\documentclass[a4paper]{article}
\usepackage[utf8x]{inputenc}
\usepackage{amssymb, amsthm, amsmath, graphicx, fancyhdr, indentfirst, newlfont, latexsym, enumitem, wasysym,  textcomp, hyperref}
\theoremstyle{remark} 
\newtheorem*{rem}{Remark}
\theoremstyle{plain} 
 \newtheorem{prop}{Proposition}
\theoremstyle{plain}
\newtheorem{teor}{Theorem}
\newtheorem*{teor*}{Theorem}
\theoremstyle{definition}
\newtheorem*{lemma}{Lemma}
\newtheorem*{cor}{Corollary}

\newcommand{\HA}{H_{\alpha}^h}
\newcommand{\HB}{H_{\beta}^h}

\newcommand{\hp}{\mathbb{H}^2}

\newcommand{\R}{\mathbb{R}}
\newcommand{\uah}{u_{\alpha}^h}

 \title{Non existence of constant mean curvature graphs on circular annuli of $\hp$}
 \author{C. Senni}
\date{March, 2011}

\begin{document}

\maketitle

\begin{abstract}

We show a non existence result for solutions of the prescribed mean curvature equation in the product manifold $\hp \times \R$, where $\hp$ is the real hyperbolic plane. More precisely we prove a-priori estimates for graphs with constant mean curvature $h  \in (0, \frac{1}{2}]$ on circular annuli of $\hp$. For $0 < h < \frac{1}{2}$ we obtain an estimate from above on any circular annulus and one from below on annuli with a small hole, the size of the hole depending on $h$. For $h = \frac{1}{2}$ we obtain both estimates for any circular annulus. All the estimates depend only on the thickness of the annulus and the value of the graph on the outer boundary.
\end{abstract}

\section*{Introduction}
In the Euclidean case it is possible to reduce the study of constant mean curvature, \textit{cmc} for short, surfaces to the two cases of zero and positive curvature, while in $\hp \, \times \, \R$ one has to distinguish at least three instances according to the value $h$ of the mean curvature: the minimal case, the case $h \, \in (0, \frac{1}{2}]$ and the case $h \, \in \, (\frac{1}{2}, +\infty)$. The role of the value $h=\frac{1}{2}$ has been outlined by Daniel on one side and Spruck for a different aspect. The first phenomenon occurring for $h = \frac{1}{2}$, discovered in \cite{BENOIT}, is the existence of a local isometry between minimal surfaces of the (Riemannian) Heisenberg group and surfaces of $\hp \times \R$ with constant mean curvature equal to $\frac{1}{2}$. The other one, described by Spruck in \cite{IGSP}, is that $h=\frac{1}{2}$ is the biggest value of the mean curvature such that the horosphere convexity of the boundary (see \cite{CIAO} and \cite{CM} for more details on this concept) is sufficient to have a solution of the constant mean curvature equation on a regular domain with any prescribed regular boundary value.\\
The study of minimal surfaces in $\hp \times \R$ was started in the 00's by the work of Nelli and Rosenberg \cite{NERO}, and was further developed by Hauswirth \cite{HAUSRIEMAN}, Meeks and Rosenberg \cite{MEEKSROSTHEORY} and \cite{MEEKSROSTABLE}, Rosenberg \cite{ROSILLY} and, more recently, Daniel \cite{BENOIT2}.  The problem of positive constant mean curvature  has been addressed by Abresch and Rosenberg in \cite{ABROHOPF} with the introduction of the generalized Hopf differential. After that Sa Earp and Toubiana in \cite{SATOU} found examples of rotational constant mean curvature surfaces, Fernandez and Mira in \cite{FEMI} proposed a construction of a Gauss Map for constant mean surfaces in $\hp \, \times \, \R$, Nelli and Rosenberg in \cite{MEME} and \cite{NEROGLO} established general theorems for these surfaces.\\
As it is well known the geometric properties of the boundary play a crucial role in the existence of solutions of the prescribed mean curvature equation. In the Euclidean case, if $\Omega$ is a regular bounded and convex domain in $\R^2$ and $\phi$ is a regular function defined on its boundary, one can always find a regular function $u$ with constant mean curvature $h \geq 0$ and prescribed boundary value $\phi$. This classical result due to Serrin \cite{SERRIN} was extended to $\hp \, \times \, \R$ and $\mathbb{S}^2 \times \R$ by Spruck \cite{IGSP} for suitable values of $h$ and boundaries satisfying a convexity condition appropriate to the Riemannian manifold considered. The same problem, but in a more general setting, was considered by Dajczer, Hinojosa and de Lira in \cite{DAHILI}, and more recently by Dajczer and de Lira in \cite{DALI}. In the hyperbolic case the geometric property of the boundary assuring existence of a solution of the prescribed constant mean curvature equation is the horosphere convexity. The prescribed mean curvature equation on a set whose boundary has an arbitrary shape in general has no solution. Already in the Euclidean case Finn \cite{FI} and Jenkins and Serrin \cite{JS}  proved that non convexity can lead to non existence.  In particular Finn in \cite{FI} considered minimal graphs on circular annuli and obtained a-priori estimates depending only on the thickness of the annulus and the value of the graph on the outer boundary.  Estimates of this kind yield non existence in that, given a circular annulus, one can assign boundary data on the inner circle violating the estimates. For these reasons it is particularly interesting to study non existence of cmc graphs on annuli of $\hp$. By circular annulus we mean a set $\Omega(a,b) = \{ z \, \in \, \hp \, : \, a \leq |z|_{\hp} \leq b \}$, where $0 <a<b$. At this stage we are not interested in regularity aspects of the problem hence we consider graphs of functions in $C^2 \Big( \mbox{int} ( \Omega(a,b) ) \Big) \, \cap \, C \Big( \Omega(a,b) \Big)$.\\
 The problem of existence and non existence of solutions of the constant mean curvature equation on annuli is of great interest for its application to the study of ends of cmc surfaces. Indeed the standard way to construct a cmc end is to consider a sequence of annuli diverging to an exterior domain, to solve the prescribed mean curvature problem on each annulus and to prove the convergence of the sequence to a solution on the exterior domain. We refer to Osserman's work \cite{OSSERMAN} and Schoen's classification \cite{BELLARIC} for an overview of the classical results in the Euclidean case. In the case of $\hp \times \R $ it is natural to start the study of cmc ends by considering values of the mean curvature $h \, \in \, (0, \frac{1}{2}]$  because in this case there are no compact closed cmc surfaces (see \cite{NEROGLO}). This problem has been addressed by Sa Earp and Toubiana \cite{SATOU} in the rotational case for $h \in (0, \frac{1}{2} ]$. Non rotational cases have been considered by Elbert, Nelli and Sa Earp \cite{NESA} for  $h = \frac{1}{2}$ and very recently by Citti and Senni in \cite{CISE} for $h \in (0, \frac{1}{2})$. \\
The main results we prove here, Theorem \ref{teor: height estimates} and Theorem \ref{teor: height estimates h=1/2}, are the a-priori estimates for functions on circular annuli whose graphs have constant mean curvature $h \, \in \, (0, \frac{1}{2}]$. Our estimates depend only on the thickness of the annulus and on the value assumed by the graph on the outer boundary, hence they do not depend on the value of the graph on the inner boundary. Moreover, all the estimates can be easily written in an explicit form. The two results take respectively into account the case $h \in (0, \frac{1}{2})$ and $h=\frac{1}{2}$. In the first case we obtain an estimate from below on any circular annulus and one from above for domains with a small hole, the size of the hole being strictly bounded by $\frac{1}{\sqrt{1-4h^2}}$. We remark that this bound blows up when $h=\frac{1}{2}$, hence in the second case we obtain both estimates for any circular annulus. In a way similar to the one proposed by Finn in \cite{FI}, our estimates depend on a Lemma which is a comparison principle for solutions of a special class of elliptic quasilinear equations, namely equations which do not explicitly depend on the value of the unknown function. This is Lemma \ref{lem: lemma finn} and is the hyperbolic analogue of Lemma 6 in \cite{FI}. Roughly speaking the Lemma states that, in order to give an estimate for a cmc $h$ graph on a circular annulus, one can use any function whose graph has cmc $h$ and is vertical in the inner boundary of the annulus. In the Euclidean case the most immediate example are the catenoids, in the hyperbolic setting we have the $\{ \HA \}_{\alpha}$ surfaces introduced by Sa Earp and Toubiana in \cite{SATOU}. \\
I would like to thank Professor Giovanna Citti and Professor Alberto Parmeggiani from the University of Bologna for the helpful suggestions about the organization of the contents.\\
The plan of the paper is the following:\\
In section one we recall some facts of hyperbolic geometry.\\
In section two we recall some properties of the rotational cmc family $\{ \HA \}_{\alpha}$.\\
In section three we prove our estimates and the corollary explicitly stating the non existence of solutions on circular annuli with appropriate boundary conditions.

\section{Hyperbolic setting}
Here we recall only the facts we are using in the paper, for a general introduction to hyperbolic geometry one can refer to \cite{ANDERSON}.\\
We consider the Poincar\'e model of the hyperbolic plane, which means that for us $\hp$ is the unit disc $ \{ z =(x,y) \, \in \R^2 : |z| < 1\}$ with the conformal metric 
\begin{align*}
d \sigma^2(z) & = \left( \dfrac{2}{1- |z|^2} \right)^2 \, \Big( dx^2 + dy^2 \Big)
\end{align*}
We denote by $\partial \hp =  \{ |z|_{\R^2}=1\}$ the asymptotic boundary of $\hp$ because it is a set at infinite distance from any point of the hyperbolic plane.
$\hp$ is a homogeneous manifold with a three dimensional group of isometries and constant sectional curvature equal to $-1$. The homogeneity allows us to choose any point as the origin.\\
 We denote by $\nabla$ the Levi Civita connection given by the metric. In this model geodesics are (suitable parametrizations of) arcs of Euclidean circles crossing orthogonally $\partial \hp$, or (suitable parametrizations of) Euclidean rays emanating from $0 \, \in \hp$. \\
To describe the Riemannian product $\hp \, \times \, \R$ we use the coordinates given by the product. Denoting $t$ as a coordinate for $\R$, the metric we are considering on $\hp \, \times \, \R$ is $ d \sigma^2 + dt^2 $. 

One can easily check that if $S \, \subset \, \hp \times \R$ is a smooth surface that is a graph on $\Omega\, \subset \, \hp$, its mean curvature $H$ can be written in terms of purely hyperbolic quantities and precisely, if $S = \{ (z, u(z)) : z \, \in \, \Omega\}$ for some $u \, \in \, C^2(\Omega)$,  we have
\begin{align*}
 2 \, H(z) = - \mathrm{div}^{\hp \, \times \, \R}(\eta(z,u(z))) = \mathrm{div}^{\hp} \, \left( \dfrac{\nabla u(z)}{\sqrt{1 + {|\nabla \, u(z)|}^2_{\hp}}} \right)
\end{align*}
where $\eta$ is the upward unit normal vector to $S$ in $\hp \, \times \, \R$. In the case of a graph, the mean curvature acts as a second order differential operator that we denote by $Q$. Precisely, if $u \, \in \, C^2(\Omega)$ we write
\begin{align*}
 2 H(z) = Q(u) =  \mathrm{div}^{\hp} \left( \dfrac{\nabla u}{\sqrt{1 + {|\nabla \, u|}^2_{\hp}}} \right)
\end{align*}

 It is well known that $Q$ is a second order operator, quasilinear and elliptic, and  uniformly elliptic whenever $|\nabla^{\hp} u|_{\hp}$ is uniformly bounded. Moreover it does not depend explicitly on $u(z)$, so reflecting the fact that each vertical translation is an isometry of $\hp \, \times \, \R$. This structure of $Q$ has the useful consequence that solutions of equations prescribing a value for $Q$ satisfy maximum and comparison principles.

\section{The $\HA$ family}
From now on $h$ will be a real number belonging to the interval $(0, \frac{1}{2} ]$. We recall the formulas and some relevant properties of the rotational cmc surfaces introduced by Sa Earp and Toubiana in \cite{SATOU}. Let $\alpha \, \in (0, +\infty)$. \\
For $0 < h \leq \frac{1}{2}$ one defines 
\begin{equation}
\rho^h(\alpha)=
\left \{
\begin{array}{cc} 
\mathrm{arccosh} \Big( \dfrac{-2 \, \alpha \, h + \sqrt{1-4h^2 + \alpha^2 }}{1 - 4h^2} \Big) & \mbox{ if } 0< h <\dfrac{1}{2} \label{eq: rho^h_alpha h < 1/2} \\
|\log(\alpha)| & \mbox{ if }  h  =\dfrac{1}{2}
\end{array}
\right.
\end{equation}
and for all $\rho > \rho^h(\alpha)$
\begin{align}
u^h_{\alpha}(\rho) & = \dfrac{-\alpha + 2h \, \mathrm{cosh}(\rho)}{\sqrt{ \mathrm{sinh}(\rho)^2 - (-\alpha + 2h \, \mathrm{cosh}(\rho))^2}} \\
\HA(\rho) & = \int_{\rho^h(\alpha)}^{\rho} \, \uah(r) \, dr
\end{align}
If $\rho$ has the meaning of the hyperbolic distance from $0 \, \in \, \hp$, the above formulas define a family of rotational surfaces in $\hp \, \times \, \R$, where by rotational we mean invariant with respect to the rotation about the line $\{ z=0 \} \, \subset \, \hp \, \times \, \R$. These surfaces are graphs defined in the complement of discs of the hyperbolic plane. The following proposition recalls some of the properties of the $\{ \HA \}$ family and its proof can be found in \cite{SATOU}, \cite{NESASATOU}.

\begin{prop}\label{prop: prop of HA}$\\$
Let $h \, \in \, (0,\frac{1}{2}]$. Then
\begin{enumerate}
\item $\forall \, \alpha > 0 $ we have
$$
Q \Big( \HA \Big) \equiv 2 h
$$
\item\label{item: raggio e parametro 1-1} $\rho^h(\alpha)$ is monotonically decreasing in the interval $(0, 2h]$  and monotonically increasing in the interval $[2h, +\infty)$. Moreover 
\begin{itemize}
\item If $0<h<\frac{1}{2}$ \label{eq: raggi piccoli h<1/2}
\begin{align*}
\rho^h \Big( (0, 2h] \Big) & = \Big[0, \mathrm{arccosh} \left( \frac{1}{\sqrt{1-4h^2}} \right)\Big) \\
\rho^h \Big( [2h, +\infty )\Big) & = [0, +\infty)
\end{align*}

\item If $h=\frac{1}{2}$
\begin{align*}
\rho \Big( (0, 2h] \Big) & = [0, +\infty) \\
\rho \Big( [2h, +\infty )\Big) & = [0, +\infty ) 
\end{align*}

\end{itemize}

\item If  $\alpha \neq 2h$, $\HA$ is zero valued and vertical on the circle of hyperbolic radius $\rho^h(\alpha)$ 
\item For $\alpha \leq 2h$ $\HA$ is nonnegative for $\rho \, \geq \, \rho^h(\alpha)$

For $\alpha \geq 2h$, $\HA$ is non-positive in a small annulus containing its boundary and positive out of this annulus. 

\end{enumerate}
\end{prop}

\begin{rem}$\\$
In view of the difference of behavior of the function $\HA$ for $\alpha < 2h$ and $\alpha > 2h$, we use the letter $\beta$ when the parameter of a surface is greater than $2h$. In other words when we write $\HB $ we tacitly assume $\beta > 2h$.
\end{rem}

\begin{figure*}[h]
\begin{center}    
\includegraphics[width=12cm]{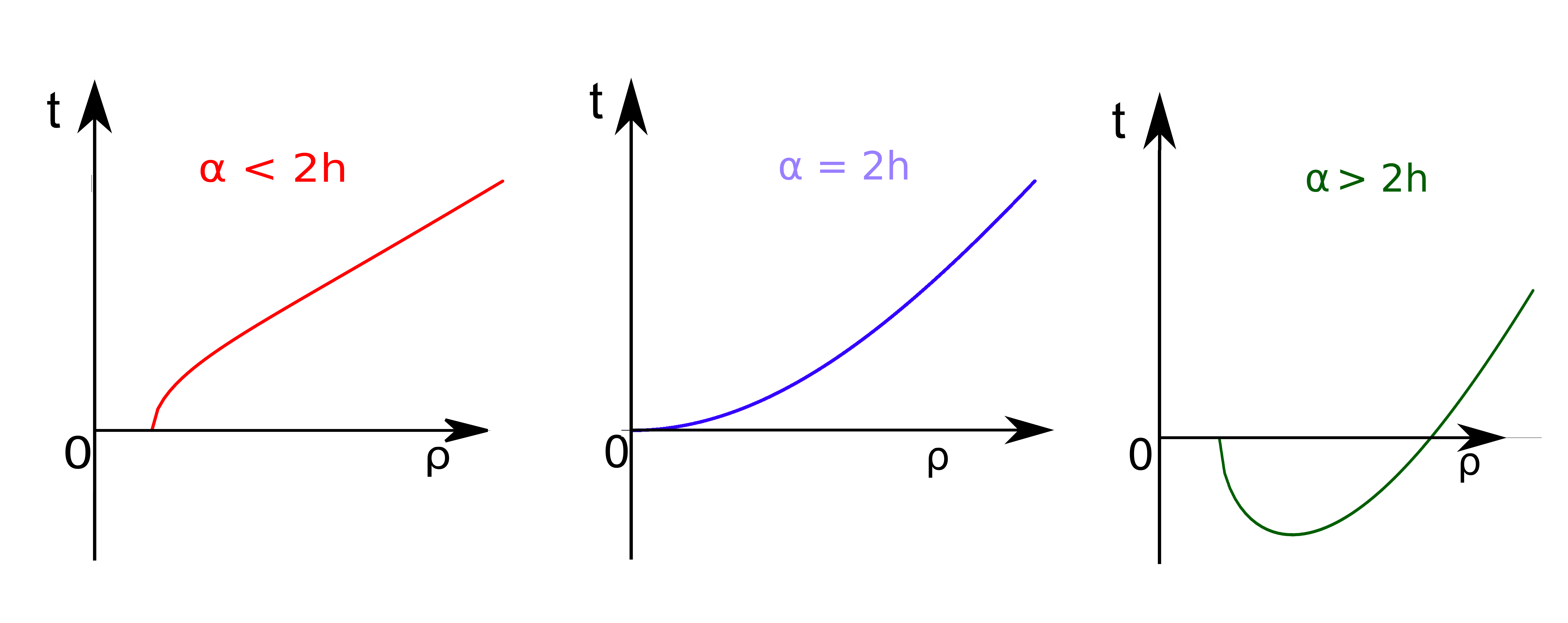}
\caption[legenda elenco figure]{The dependence of $\HA$ on the parameter $\alpha$}\label{fig: sovrapposte}
\end{center}
\end{figure*}

Figure \ref{fig: sovrapposte} shows the dependence of the shape of the generating curves on the parameter $\alpha$.

\section{A-priori estimates and non existence results}
In this section we prove the a-priori estimates for graphs with constant mean curvature $h \in \left( 0, \frac{1}{2} \right]$ on circular annuli of $\hp$. As we have mentioned, our estimates follow from a Lemma which can be proved under very general hypotheses. Before introducing the Lemma let's define some notation. Consider $\Omega \, \subset \, \hp$ a compact annulus. Let $\partial \Omega = \gamma_1 \cup \gamma_2$, where $\gamma_i$ is a Jordan curve for $i=1,2$ and assume $\gamma_1$ is contained in the compact set bounded by $\gamma_2$, note that $\gamma_1$ and $\gamma_2$ can have non empty intersection. We are not going to require that the boundary of $\Omega$ is smooth, we only require the possibility of defining a derivative of functions defined in $\mbox{int} ( \Omega )$ near $\gamma_1 \subset \partial \Omega$. To do this we suppose that each $z \, \in \, \gamma_1 $ is the end-point of a geodesic contained in $\mbox{int} ( \Omega )$. More precisely, for any given $ z \, \in \, \gamma_1$ we require the existence of a geodesic $\gamma : [0,a] \to \hp$ such that $\gamma \Big( [0,a) \Big) \, \subset \, \mbox{int} (\Omega )$ and $\gamma(a)=z$, for some $a>0$. 
We denote by $\frac{ D}{d s}$ the covariant derivative associated to $\gamma$, provided  that $\gamma$ ends on $\gamma_1$. 

\begin{rem}$\\$
Circular annuli satisfy the conditions just described. Indeed, chosen a point on any of the two circles bounding the annulus one can get to that point by  following a geodesic ray.
\end{rem}
%
%
%
%

\begin{lemma}\label{lem: lemma finn}$\\$
Assume $\Omega$ is a domain satisfying the hypotheses just described. Consider $\phi_0,\phi : \Omega \to \R$  two functions such that:
\begin{enumerate}
 \item $\phi_0 \, \in \, C^2 \Big( \mbox{int} (\Omega ), \R \Big) $ 
    \subitem $\phi_0$ has a limit, finite or infinite, on $\gamma_1$
    \subitem $\lim_{z \to \gamma_1} \frac{ D \, \phi_0}{d s} (z) = + \infty$ along any given geodesic ending on $\gamma_1$
\item $\phi \, \in \, C^2 \Big( \mbox{int} (\Omega ) , \R \Big) \, \cap  \, C \Big( \Omega \Big)$
\item $\liminf_{z \rightarrow \gamma_2} (\phi_0 - \phi) \geq 0$
\item $F(\phi_0) \leq F(\phi)$ in $\mbox{int} (\Omega )$, where $F(u)$ is a quasilinear and elliptic second order operator with coefficients not explicitly depending  on $u$.
\end{enumerate}
Then 
\begin{align*}
 \liminf_{z \rightarrow \gamma_{1}} (\phi_0 - \phi) (z) \geq 0
\end{align*}
\end{lemma}

\begin{rem}$\\$
\begin{itemize}
\item The proof of the Euclidean version of this Lemma can be found in \cite{FI}. This proof is based on topological properties of $\R^2$ and structure properties of the mean curvature operator hence it can be used also in the hyperbolic case.
 \item This Lemma can be proved for more general domains. It would be enough to require that the $\gamma_i$ are closed sets and that the derivative of $\phi_0$ blows up only in the interior points of $\gamma_1$.  However the hypotheses we made are general enough to obtain interesting results.
\end{itemize}
\end{rem}

We now show the a-priori estimates for cmc graphs on circular annuli. Our estimates depend only on the thickness of the annulus and on the boundary value of the graph on the outer boundary. To prove the estimates we apply the Lemma. Roughly speaking, what we do is to associate to an annulus the  $\HB$ and the $\HA$ defined on the complement of the disc bounded by the annulus, and then  find suitable vertical translations giving the estimates.\\
As we have mentioned (see Proposition \ref{prop: prop of HA}) the elements of $\{ \HA \}_{\alpha}$ have a different behavior depending on whether $h \in (0, \frac{1}{2})$ or $h = \frac{1}{2}$.  Hence it is very natural to distinguish two cases according to the value of the mean curvature: in the $h \, \in (0, \frac{1}{2} )$ case we obtain an estimate from above for any annulus and one from below only for annuli with a small hole. In the case $h=\frac{1}{2}$ we obtain both estimates for any circular annulus. This is because for $h \, \in (0, \frac{1}{2} )$ the $\HA$ for $0 < \alpha < 2h$ are defined on complements of small discs (recall the dependence of the radius on the parameter \eqref{eq: raggi piccoli h<1/2}), while for $h = \frac{1}{2}$ they are defined on the complement of any disc.\\
 Given $0<a<b$ we denote a circular annulus by $\Omega(a,b) = \{ z \, \in \, \hp : a \leq |z|_{\hp} \leq b \}$.
\begin{teor}\label{teor: height estimates}$\\$
 Let be $0 < h < \frac{1}{2}$ and consider $u  \, \in \, C^2 \Big( \mathrm{int} ( \Omega(a,b) ) \Big) \, \cap \, C \Big( \Omega (a,b) \Big)$ such that
$$
Q(u) = 2\, h
$$
Then for all $z \, \in \, \Omega(a,b)$   we have
\begin{align}
 u(z)  \leq H_{\beta}^h(z) -   H_{\beta}^h(b) + M \label{eq: finn da sopra} 
\end{align}
where 
$$
\beta = (\rho^h)^{-1}(a)  \, \in (2h, +\infty) \quad \mbox{ and } \quad M = \max_{ \{ |z|_{\hp} = b \} } u
$$
If moreover $a < \frac{1}{\sqrt{1-4h^2}}$, we have
\begin{align}
 \HA(z) - \HA(b) + m \leq u(z)   \label{eq: finn da sotto} 
\end{align}
where
$$
\alpha = (\rho^h)^{-1}(a) \, \in \, (0,2h) \quad \mbox{ and } \quad m=\min_{ \{ |z|_{\hp} = b \}} u
$$

\proof $\\$
Let's prove the first inequality. Recall that, by item \ref{item: raggio e parametro 1-1} of Proposition \ref{prop: prop of HA}, for any given radius $\rho^*$ there is a $\beta > 2h$ such that $\HB$ is defined in the complement of the circle of radius $\rho^*$. By equation \eqref{eq: rho^h_alpha h < 1/2} this $\beta(\rho^*)$ is given by
$$
\beta(\rho^*) = 2h \, \cosh(\rho^*) + \sqrt{\cosh(\rho^*)^2 -1}
$$

Since the differentiability of $u$ is not assumed on the boundary, we show that the claim holds on circular annuli $\Omega(a^*,b)$ for each $ a^* \, \in (a,b)$. The claim will follow by continuity of $u$ up to the boundary. Let's start by proving inequality \eqref{eq: finn da sopra}. For all $ z \, \in \, \Omega(a^*, b)$ we define
\begin{align*}
\beta & = \beta(a^*) \\
 \phi_0(|z|_{\hp}) & =  H^h_{ \beta(a^*)}(|z|_{\hp}) - H^h_{ \beta(a^*)}(b) + M
\end{align*}
$\HB$ being rotational, we obtain
$$
 u_{\{ |z|_{\hp} = b\}} \leq {\phi_0}_{| \{ |z|_{\hp} = b\}}  = M
$$
To have the same inequality on  $\{ |z|_{\hp} = a^* \}$ we observe that $H_{\beta(a^*)}^h$ is vertical on $\{ |z|_{\hp} = a^*\}$ and negative nearby this boundary, hence 
$$
\lim_{p \rightarrow \lbrace |z|_{\hp} = a^*\rbrace}  \dfrac{ D \, \HB}{d s} (p) = + \infty
$$
Moreover by the hypotheses we have
$$
Q(\phi_0) = 2 h = Q(u)
$$
and thus, $\phi_0$ being radial, the Lemma yields
$$
 u_{| \{ |z|_{\hp} = a^*\}} \leq \phi_0(a^*)
$$
Applying the standard maximum principle we obtain the inequality
$$
u(z) \leq \phi_0(|z|_{\hp}) \quad \mbox{ on } \; \Omega(a^*,b)
$$
Now suppose $a < \frac{1}{\sqrt{1-4h^2}}$. As in the preceding case recall that, by item \ref{item: raggio e parametro 1-1} of Proposition \ref{prop: prop of HA}, there exists $\rho^* \in \left( a,  \dfrac{1}{\sqrt{1-4h^2}} \right)$ such that there is a unique $\alpha < 2h$ so that $\HA$ is zero on the circle of radius $\rho^*$. By equation \eqref{eq: rho^h_alpha h < 1/2} this $\alpha(\rho^*)$ is given by
$$
\alpha(\rho^*) = 2h \, \cosh(\rho^*) - \sqrt{\cosh(\rho^*)^2 -1}
$$

To prove inequality \eqref{eq: finn da sotto} we define $\forall \,  a^* \, \in (\rho^*, b)$
$$
\phi_1(|z|_{\hp}) =  H_{\alpha(a^*)}^h(|z|_{\hp}) - H_{\alpha(a^*)}^h(b) + m 
$$
and repeat the argument of the preceding case. Applying the Lemma to $-\phi_1$ we obtain $\forall \, a^* \, \in \, (a, b)$
$$
\phi_1(a^*) \leq u(z)_{|\{ |z|_{\hp} = a^*\}} 
$$
which completes the proof.
\endproof
\end{teor}

\begin{rem}$\\$
 An interesting feature of the above result is the difference with the Euclidean minimal case. Indeed in that case one can use catenoids to establish both estimates, the one from below and the one from above (see \cite{FI}). This is because if the Euclidean mean curvature operator applied to $u$ gives zero, the same holds for $-u$. This is not anymore true when the mean curvature has a sign. For example in our case if  $u$ is defined on a subset of $\hp$ and $Q(u) = 2 h$, we have $Q(-u) = - 2h$ and hence we cannot use the same surface to obtain the estimate from below and the estimate from above. To cope with this problem we use the fact that the elements of the family $\{ \HA \}_{\alpha}$ approach their boundary with negative singular normal derivative when $ \alpha \in (0,2h )$ and positive singular normal derivative when $ \alpha \in (2h,+ \infty )$.
\end{rem}

Even if we can give an a-priori bounding box only for circular annuli with a small hole, we have a non existence result for any annulus. This is the content of the next Corollary

\begin{cor}$\\$
 Let $\Omega(a,b)$ be a circular annulus in $\hp$. Then $\forall \, c>0$ and $\forall  \, \varepsilon > 0$ the following Dirichlet problem has no solution in $C^2 \Big(\text{int} (\Omega(a,b) ) \Big) \cap C \Big( \Omega(a,b) \Big)$.

$$
\left \{ \begin{array}{cll}
Q(u) & = 2h & \mbox{ in } \Omega(a,b) \\
u & =   \HB(b) + c + \varepsilon   &  \mbox{ in } \{ |z| = a\}  \\
u & =   c    & \mbox{ in } \{ |z| = b\} 
         \end{array}
\right.
$$

\proof $\\$
In this case $\forall \, z \, \in \, \{ |p|_{\hp} = a \}$  Theorem \ref{teor: height estimates} gives
$$
u(z)  \leq    H_{\beta}^h(b) + c
$$
hence the boundary data on $\{ |z| = a\}$ cannot be achieved.
\endproof
\end{cor}

We end the section with the a-priori estimates for the $h=\frac{1}{2}$ case. The proof is as in Theorem 1

\begin{teor}\label{teor: height estimates h=1/2}$\\$
 Let  $ h = \frac{1}{2}$ and consider $u  \, \in \, C^2 \Big( \mathrm{int} (\Omega(a,b) ) \Big) \cap C \Big( \Omega (a,b) \Big)$ such that
$$
Q(u) = 2\, h = 1
$$
Then for all $ z \, \in \, \Omega(a,b)$   we have
\begin{align}
H^{h}_{\alpha}(z) - H^h_{\alpha}(b) + m \leq  u(z)  \leq H^h_{\beta}(z) -   H^h_{\beta}(b) + M \label{eq: finn h=1/2} 
\end{align}
where 
\begin{align*}
\alpha & = ( \rho^h )^{-1}(a) \, \in \, (0,1) \quad  & \quad \beta  & = (\rho^h)^{-1}(a)  \, \in (1, +\infty)
\intertext{and}
m & =\min_{ \{ |z|_{\hp} = b \}} u  \quad  & \quad M  & = \max_{ \{ |z|_{\hp} = b \} } u
\end{align*}

\proof $\\$
The proof can be done in the very same way as in the Theorem $1$. We only observe that since
 $$
\rho^{\frac{1}{2}}(\alpha) =
 \left\{
 \begin{array}{cccc}
& -\log(\alpha)  & \mbox{ for  } & 0 < \alpha \leq 1 \\
& \log(\alpha)  & \mbox{ for  }  &  \alpha \geq 1 \\
 \end{array}
 \right.
 $$ 

we can choose an $\HA$ defined on the complement of any disc, hence the estimate from below holds on any circular domain.
\endproof
\end{teor}

\begin{figure*}[ht]
\begin{center}    
\includegraphics[width=9cm]{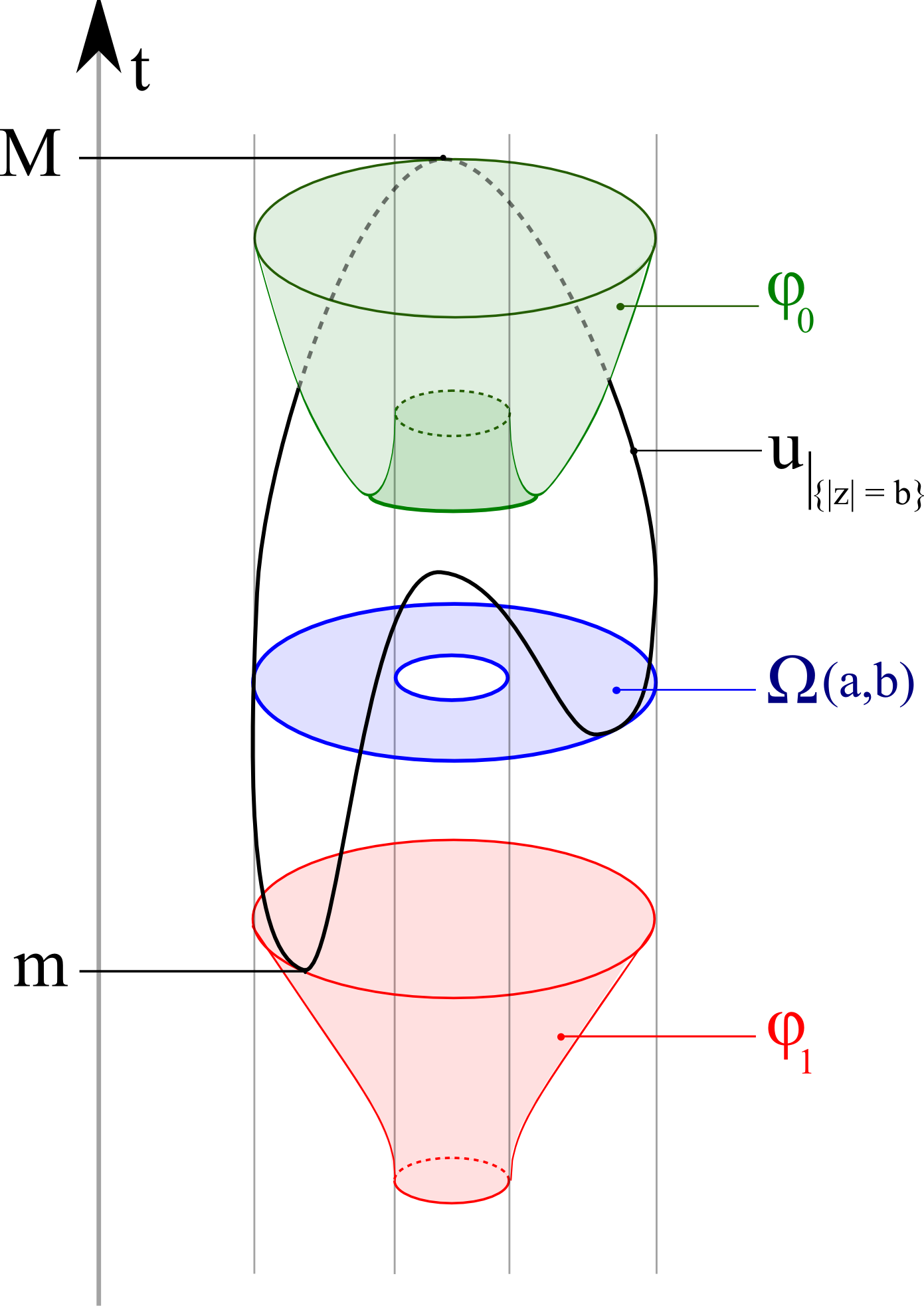}
\caption{The bounding box of Theorem \ref{teor: height estimates h=1/2}}\label{fig: finn}
\end{center}
\end{figure*}

Figure \ref{fig: finn} shows the bounding box we have just built.

\newpage
\bibliography{biblio}

\providecommand{\bysame}{\leavevmode\hbox to3em{\hrulefill}\thinspace}
\providecommand{\MR}{\relax\ifhmode\unskip\space\fi MR }
\providecommand{\MRhref}[2]{%
  \href{http://www.ams.org/mathscinet-getitem?mr=#1}{#2}
}
\providecommand{\href}[2]{#2}
\begin{thebibliography}{10}

\bibitem{ANDERSON}
J.~W. Anderson, \emph{Hyperbolic geometry}, second ed., Springer, 2005.

\bibitem{CM}
E.~Cabezas-Rivas and V.~Miquel, \emph{Volume preserving mean curvature flow in
  the hyperbolic space}, Indiana Univ. Math. J. \textbf{56} (2007), 2061--2086.

\bibitem{CISE}
G.~Citti and C.~Senni, \emph{Constant mean curvature graphs on exterior domains
  of the hyperbolic plane}, preprint (2010).

\bibitem{CIAO}
R.~J. Currier, \emph{On hypersurfaces of hyperbolic space infinitesimally
  supported by horospheres}, Trans. American Math. Soc. (1989), no.~313,
  419--431.

\bibitem{DALI}
M.~Dajczer and J.~H. de~Lira, \emph{Killing graphs with prescribed mean
  curvature and riemannian submersions}, Ann. Inst. H. Poincaré Anal. Non
  Lin\'eaire \textbf{26} (2009), no.~3, 763–775.

\bibitem{DAHILI}
M.~Dajczer, P.~A. Hinojosa, and J.~H. de~Lira, \emph{Killing graphs with
  prescribed mean curvature}, Calc. Var. Partial Differential Equations
  \textbf{33} (2008), no.~2, 231–248.

\bibitem{BENOIT}
B.~Daniel, \emph{Isometric immersions into 3-dimensional homogeneous
  manifolds}, Comment. Math. Helv. \textbf{82} (2007), no.~1.

\bibitem{BENOIT2}
\bysame, \emph{Isometric immersions into {$\mathbb{S}^n \times \mathbb{R}$} and
  {$\mathbb{ H} ^n \times \mathbb{R}$} and applications to minimal surfaces},
  Trans. Amer. Math. Soc. (2009), no.~12, 6255--6282.

\bibitem{FI}
R.~Finn, \emph{Remarks relevant to minimal surfaces, and to surfaces of
  prescribed mean curvature}, J. d'Anal. Math. \textbf{14} (1965), no.~1,
  139--160.

\bibitem{HAUSRIEMAN}
L.~Hauswirth, \emph{Minimal surfaces of riemann type in three-dimensional
  product manifolds}, Pacific J. Math. \textbf{224} (2006), no.~1, 91--117.

\bibitem{JS}
H.~Jenkins and J.~Serrin, \emph{The dirichlet problem for the minimal surface
  equation in higher dimensions}, J. Reine Angew. Math. \textbf{229} (1968),
  170--187.

\bibitem{MEEKSROSTABLE}
W.~H.~III Meeks and H.~Rosenberg, \emph{Stable minimal surfaces in {$M\times
  \mathbb{R}$}}, J. Differential Geom. \textbf{68} (2004), no.~3.

\bibitem{FEMI}
P.~Mira and I.~Fernandez, \emph{Harmonic maps and constant mean curvature
  surfaces in {$\Bbb H^2\times\Bbb R$}}, Amer. J. Math. \textbf{129} (2007),
  no.~4, 1145--1181.

\bibitem{MEME}
S.~Montaldo and F.~Mercuri, \emph{A weierstrass representation formula for
  minimal surfaces in {$\mathbb H_3$ and $\mathbb H^2\times\mathbb R$}}, Acta
  Math. Sin. (Engl. Ser.) \textbf{22} (2006), no.~6, 1603--1612.

\bibitem{NEROGLO}
B.~Nelli and H.~Rosenberg, \emph{Global properties of constant mean curvature
  surfaces in {$\mathbb H^2 \times \mathbb R$}.}, Pacific J. Math. \textbf{226}
  (2006), no.~1, 137--152.

\bibitem{NESA}
B.~Nelli and R.~{Sa Earp}, \emph{Vertical ends of constant mean curvature {$H =
  \frac{1}{2}$} in {$ \mathbb{H}^2 \times \mathbb{R}$}}, Preprint (2007).

\bibitem{NESASATOU}
B.~Nelli, R.~{Sa Earp}, W.~Santos, and E.~Toubiana, \emph{Uniqueness of
  {$H$}-surfaces in {$ \mathbb{H}^2 \times \mathbb{R}$} , {$\vert H\vert
  \le1/2$}, with boundary one or two parallel horizontal circles.}, Annals of
  Global Analysis and Geometry \textbf{33} (2008), no.~4, 307--321.

\bibitem{OSSERMAN}
R.~Osserman, \emph{Global properties of minimal surfaces in {$E^3$} and
  {$E^n$}}, Annals of Mathematics \textbf{80} (1964), no.~2, 340--364.

\bibitem{ROSILLY}
H.~Rosenberg, \emph{Minimal surfaces in {${\mathbb{M}}^2\times \mathbb{R}$}},
  Illinois J. Math. \textbf{46} (2002), 1177--1195.

\bibitem{ABROHOPF}
H.~Rosenberg and U.~Abresch, \emph{A hopf differential for constant mean
  curvature surfaces in {$S^2 \times R$ and $H^2 \times R$}}, Acta Math.
  \textbf{193} (2004), no.~2, 141--174.

\bibitem{MEEKSROSTHEORY}
H.~Rosenberg and W.~H.~III Meeks, \emph{The theory of minimal surfaces in
  {$M\times \mathbb{R}$}}, Comment. Math. Helv. \textbf{80} (2005), no.~4,
  811--858.

\bibitem{NERO}
H.~Rosenberg and B.~Nelli, \emph{Minimal surfaces in {${\mathbb H}^2 \times
  \mathbb{R}$}}, Bull. Braz. Math. Soc. (N.S.) \textbf{33} (2002), no.~2,
  263--292.

\bibitem{SATOU}
R.~{Sa Earp} and E.~Toubiana, \emph{Screw motion surfaces in {$ \mathbb{H}^2
  \times \mathbb{R} $} and {$ \mathbb{S}^2 \times \mathbb{R} $}}, Illinois
  Journal of Mathematics \textbf{49} (2005), no.~4, 1323--1362.

\bibitem{BELLARIC}
R.~Schoen, \emph{Uniqueness, symmetry, and embeddedness of minimal surfaces},
  J. Differential Geom. \textbf{18} (1983), no.~4, 791--809.

\bibitem{SERRIN}
J.~Serrin, \emph{A priori estimates for solutions of the minimal surface
  equation.}, Arch. Rational Mech. Anal. \textbf{4} (1963), 376--383.

\bibitem{IGSP}
J.~Spruck, \emph{Interior gradient estimates and existence theorems for
  constant mean curvature graphs in {$ M^n \times \mathbb{R} $}}, Pure Appl.
  Math. Q. (2007), no.~3, 785--800.

\end{thebibliography}
\bibliographystyle{amsplain}

\end{document}